\documentclass[10pt]{article}

\usepackage{latexsym}
\usepackage{amssymb}
\usepackage{amsmath}

\newtheorem{Theorem}{Theorem}[section]

\newtheorem{thm}[Theorem]{Theorem}

\newtheorem{proposition}[Theorem]{Proposition}
\newtheorem{prop}[Theorem]{Proposition}

\newtheorem{lemma}[Theorem]{Lemma}

\newtheorem{cor}[Theorem]{Corollary}

\newtheorem{definition}[Theorem]{Definition}

\newcommand{\Imp}{\mbox{$\Longrightarrow$}}
\newcommand{\Iff}{\mbox{$\Longleftrightarrow$}}

\def\proof{{\bf Proof.}\ }

\def\ie{\emph{i.e.}}
\def\pes{\emph{e.g.}}
\def\Pes{\emph{E.g.}}
\def\A{{\mathcal{A}}}

\def\E{{\mathcal E}}
\def\F{{\mathcal F}}

\def\G{\mathcal{G}}

\def\P{{\mathcal P}}
\def\S{{\mathcal S}}

\def\U{{\mathcal U}}

\def\C{{\mathcal C}}

\def\N{{\mathbb N}}

	\def\B{\mathcal{B}}

\def\fg{{\varphi}}

\def\og{{\omega}}

\def\sg{{\sigma}}

\def\cof{\mbox{\rm cof}\;}
\def\cov{\mbox{\bf cov}}

\def\min{\mbox{\rm min}\;}
\def\max{\mbox{\rm max}\;}

\def\+#1{\vec{#1}}

\def\ck{{\mathfrak{c}}}

\def\dk{{\mathfrak{d}}}

\def\Sk{{\mathfrak{S}}}

\def\*{\times}

\def\0{\emptyset}
\def\/{\setminus}
\def\_{\overline}

\def\<{\prec}
\def\Equ{\equiv_{\U}}

\def\ult#1#2{^{#1}_{\; #2}}

\def\incl{\subseteq}

\def\la{\langle}
\def\ra{\rangle}

\def\zfc{\textsf{ZFC}}
\def\ZFC{\textsf{ZFC}}
\def\CH{\textsf{CH}}
\def\MA{\textsf{MA}}

\def\qed{${}$\hfill $\Box$}

\def\etc{\emph{etc.}}

\begin{document}

\title{Quasi-selective and
\\
weakly Ramsey ultrafilters}

\author{Marco Forti\\
{\small Dipart. di Matematica Applicata ``U. Dini'',
Universit\`{a} di Pisa, Italy.}\\
{\small \tt{forti@dma.unipi.it}}}
\date{}
\maketitle


\begin{abstract}
\emph{Selective}  ultrafilters are characterized by many equivalent 
properties, in particular the Ramsey property that
every finite colouring of $[\N]^{2}$ has a \emph{homogeneous} set $U\in\U$, 
and the equivalent property that
 every function is \emph{nondecreasing} on some $U\in\U$.
Natural weakenings of these properties led to the inequivalent notions 
of  \emph{weakly 
Ramsey}   and of \emph{quasi-selective} ultrafilter, introduced and 
studied in  \cite{bl74} and \cite{bdf}, respectively. Call
$\,\U$  \emph{weakly Ramsey}  
    if for
    every finite colouring of $[\N]^{2}$ there is $U\in\U\,$
    s.t. $[U]^{2}$ 
    has \emph{only two colours}, and call
     $\,\U$  $f$-\emph{quasi-selective} if every function
$g\le f$ is \emph{nondecreasing} on some $U\in\U$.
 (So the \emph{quasi-selective} ultrafilters of \cite{bdf} are here
 $id$-quasi selective.)
In this paper we consider the relations between various natural cuts 
of the ultrapowers of $\N$ modulo  weakly Ramsey 
and $f$-quasi-selective ultrafilters. In particular we characterize 
those weakly Ramsey ultrafilters that are isomorphic to a 
quasi-selective ultrafilter.

\end{abstract}

\maketitle

\section*{Introduction}

Special classes of ultrafilters over $\N$ have been
introduced and variously applied in the literature,
starting from the pioneering work by G. Choquet \cite{co1,co2}
in the sixties (see \pes\ \cite{bo}). 
Particular attention received the class of \emph{selective} (also called 
\emph{Ramsey}, or in French \emph{absolute}) ultrafilters.
It is well known that the ultrafilter $\U$ is 
selective  if and only if 
every finite colouring of $[\N]^{2}$ has a homogeneous  set $U\in\U$
(\ie\ $[U]^{2}$ is momochromatic), 
or equivalently 
if and only if every function $f:\N\to\N$ is nondecreasing on some 
$U\in\U$.

Allowing  sets $U$ such that $[U]^{2}$ is \emph{dichromatic} in the first characterization led to the notion
of  \emph{weakly 
Ramsey} ultrafilter over $\N$, 
introduced and 
studied in  \cite{bl74} (see also \cite{ros}).
On the other hand, restricting the second characterization to functions 
\emph{smaller than the identity} defines the \emph{quasi-selective} ultrafilters over $\N$, 
introduced and 
studied in  \cite{bdf}. 
Quasi-selective ultrafilters have
 independent interest, because they are necessary in modelling the 
``Euclidean numerosities'' of point sets considered in 
\cite{bdf}, as well as in providing the so called ``fine densities'' 
of sets of natural numbers in \cite{dn}.

In this paper we make a comparative study of \emph{weakly Ramsey}  
and \emph{f-quasi-selective} 
ultrafilters, the latter class being the natural parametric generalization 
of 
quasi-selective ultrafilters, where a function $f:\N\to\N$ replaces
the identity in the original definition of \cite{bdf}.

It is worth mentioning that, on the one hand,  selective 
ultrafilters
are simultaneously weakly Ramsey and quasi-selective, while in turn 
both these classes are P-points.
 On the other hand these classes are distinct, provided that there exist
 a selective and a non-selective 
quasi-selective ultrafilter.
 The existence of these ultrafilters is not provable in \ZFC, but 
follows from mild set theoretical hypotheses, \pes\ 
the Continuuum Hypothesis \CH, or Martin's Axiom \MA. The 
study of weak sufficient conditions for the existence of all the
various kinds of these ultrafilters seems to be an interesting 
field of set theoretic research, very little explored up to now.

\medskip
The paper is organized as follows. In Section \ref{fqsu} we introduce 
the class of $f$-quasi-selective ultrafilters on $\N$, and we study their 
properties, generalizing some results of \cite{bdf}.
In section \ref{wru} we study the weakly Ramsey ultrafilters 
introduced in \cite{bl74}, and we give a complete classification in 
terms of the mutual ordering of three natural cuts of the 
corresponding ultrapowers of $\N$. We thus specify also 
the respective  properties of ``quasi-selectivity''.
Final remarks and open questions may be found in the concluding 
section \ref{froq}.

\smallskip
In general, we refer to \cite{CK} and \cite{bl} for definitions and basic facts
concerning ultrafilters and ultrapowers.

\smallskip
The author is grateful to  Mauro Di Nasso for many useful discussions, 
and to Andreas Blass for some basic suggestions.

\section{$f$-quasi-selective  ultrafilters}\label{fqsu}

Throughout this paper $\U$ is a nonprincipal ultrafilter  on $\N$, 
and all functions are $\N\to\N$, unless  different mention is 
made explicitly. 
Recall that two functions $f,g$ are $\U$-equivalent 
(written $g\Equ f$) if there exists $U\in\U$ such that $f(u)=g(u)$ 
for all $u\in U$.
In general we say that a function $f$ is increasing, unbounded, 
one-to-one, \etc, \emph{modulo} $\U$ if there exists $U\in\U$ such that the 
restriction of $f$ to $U$ is increasing, unbounded, 
one-to-one, \etc.

\begin{definition}\emph{
    Let $\U$ be a nonprincipal ultrafilter  on $\N$, and let 
    $f:\N\to\N$ be unbounded modulo $\U$. Then}
    \begin{itemize}
        \item \emph{~$\U$ is $f$-\emph{quasi-selective} (shortly $f$-QS) if,
 for all $g:\N\to\N$,} 
 $$\exists U\in\U\,\forall x\in U\, g(x)\le f(x)\ 
 \ \Imp\ \ g\ \mbox{\emph{nondecreasing}} 
 \!\!\!\mod \U.$$ 
    
        \item  \emph{~$\U$ is \emph{quasi-selective} (shortly QS) if it is 
$id$-QS, where $id:\N\to\N$ is the identity. }
    
        \item \emph{~$\U$ is \emph{properly quasi-selective} (shortly PQS) if it 
is $f$-QS for some, but not for all functions $f$. }

 \item \emph{~$\U$ is \emph{strongly quasi-selective} (shortly SQS) if it 
is $f$-QS for some $1$-$1$ function $f$. $\U$ is \emph{weakly quasi-selective} 
(shortly WQS) if it 
is PQS, but not SQS.}
    \end{itemize}

\end{definition}

Clearly the ultrafilter $\U$ is selective if and only if it is 
$f$-QS for all $f$.

Recall that the ultrafilter $f\U$ is defined by $f\U=\{V\mid 
f^{-1}[V]\in\U\}$.

Useful relations between QS ultrafilters and generic 
    $f$-QS ultrafilters are given in the following proposition:

\begin{prop}\label{fQS} ${}$
   \begin{enumerate}
    \item   If $\,\U$ is $f$-QS, then 
$f\U$ is QS.

    \item If $f$ is increasing modulo $\U$, then
     $\U\ is\ (g\circ f)\mbox{-QS}$  if and only if   
$f\U$ is $g$-QS;
in particular\ $\U$ is $f$-QS  if and only if
$f\U$ is QS.

\end{enumerate}

\end{prop}

\proof ${}$

\medskip
\noindent
1. Let $\U$ be $f$-QS, with $f$ nondecreasing on $U\in\U$. Assume that
$h(x)\le x$ for $x\in f[V], V\in\U$, so that $h\circ f\le  f$ on $U\cap V$.
Then both $f$ and $h\circ f$ are nondecreasing on $U\cap V$.
Suppose by contradiction that there exist $x,y\in U\cap V$ such that 
$f(x)<f(y)$, but $h(f(x))>h(f(y))$: the first inequality implies 
$x<y$, whereas the second implies $x>y$, contradiction. Therefore $h$ 
is nondecreasing on $f[U\cap V]\in f\U$.

\medskip
\noindent
2. Pick $U\in\U$ such that, for all $x,y\in U $, $x<y\ \Iff\ f(x)<f(y)$. 
Then, for every function $h$,
$$ \forall x,y\in U\ (\,x<y \ \Imp \ h(x) \le h(y)\,)$$
is equivalent to 
$$ \forall z,w\in f[U]\ (\,z<w \ \Imp \ h(f^{-1}(z)) \le h(f^{-1}(w))\,).$$
Moreover
$$\forall x\in U.\; h(x)<g(f(x)) \ \Iff\ \forall z\in f[U].\; 
h(f^{-1}(z))<g(z).$$
So, if $f\U$ is $g$-QS and $h<g\circ f$ on $U$, then $h\circ 
f^{-1} <g $ on $f[U]$, and hence $h\circ 
f^{-1}$ is nondecreasing on $f[U]$, which in turn is equivalent to $h$ 
nondecreasing on $U$.

Similarly, if $\U$ is $(g\circ f)$-QS and $h<g$ on $f[U]$, then $h\circ 
f <g\circ f $ on $U$, so $h\circ f$ is nondecreasing on $U$, and $h=h\circ 
f\circ f^{-1}$ is nondecreasing on $f[U]$.

The last assertion is the case  $g=id$.

\qed

It is proved in \cite{bdf} that, when $\U$ is QS, every function is 
$\U$-equivalent either to a constant, or to an  
 ``interval-to-one'' function, \ie\ a function $g$  such that, 
for all $n$,
$g^{-1}(n)$ is a (finite, possibly empty) \emph{interval} of $\N$. 
A weaker property, still sufficient to imply P-pointness, holds for 
all PQS ultrafilters, namely:

\begin{prop}\label{part}
    Let $\U$ be a PQS ultrafilter and let $\la X_{n}\mid n\in\N\ra$
    be a partition of 
    $\N$ such that no part $X_{n}$ is in $\U$. Then there exists
    an \emph{interval partition} $\la Y_{m}\mid m\in\N\ra$ and a set 
    $U\in\U$ such that
   $$\forall n \,\exists m\ X_{n}\cap U\incl Y_{m}. $$
    In particular every 
    function is either constant or
 ``finite-to-one'' modulo $\U$. 
Hence all PQS ultrafilters are nonselective P-points.

 \end{prop}   

 \proof
    Let $f$ be a nondecreasing unbounded function such that 
    $\U$ is $f$-QS. Define the function $g$ by 
    $$g(x)= f(\min X_{n})=\min f(X_{n})\ \ \mbox{for all}\ \ 
    x\in X_{n}.$$
      Then $g\le f$, so there exists a nondecreasing function $h$ 
      that is equal to $g$ on some set $U\in\U$. 
    The partition
    $\la Y_{m}=h^{-1}(m)\mid m\in\N\ra$ is 
    an interval partition that satisfies the wanted condition, 
    because $h$ is constant on $X_{n}\cap U$.
    
    \qed
    
    \medskip
    Remark that if $f$ is one-to-one, then each nonempty $Y_{m}\cap U$ is 
    equal to one $X_{n}\cap U$. In particular, modulo a SQS 
    ultrafilter, every non-constant function is interval-to-one. 
    
    \medskip

\medskip   
Recall that the ultrafilter $\U$  is \emph{rapid} if for
 every increasing function $g$ there exists 
 $U=\{u_{1}<u_{2}<\ldots<u_{n}<\ldots\}\in\U$ such that 
 $u_{n}>g(n)$. If moreover 
    $\U$ is a P-point, then $\U$ is rapid if and only if the 
    functions that are $1$-to-$1$ modulo $\U$ are coinitial in the 
    nonstandard part of the ultrapower $\N\ult{\N}{\U}$ (see \pes\ 
    \cite{blm}). It is well known that the existence of nonselective rapid 
    P-points is consistent, see \pes\ \cite{BC}. However these cannot be PQS 
    ultrafilters, since we have
 
\begin{prop}\label{rap}
     Let $\U$ be $f$-QS: then $\U$ is rapid if and only if
 it is selective. 

\end{prop}

\proof ${}$ Every selective ultrafilter is rapid, so we have to prove 
the `only if' part.
Let $\U$ be $f$-QS and let $\P=\{[p_{n},p_{n+1})\mid n\in\N\,\}$ be an interval 
partition of $\N$. By possibly unifying some intervals, we may
assume w.l.o.g. that $f(p_{n})>n$. By rapidity, there is a set 
$U=\{u_{1}<u_{2}<\ldots<u_{n}<\ldots\}\in\U$ such that $u_{n}>p_{n}$. 
  Define the function $g$ by $$g(x)=|\{m\le n\mid x\le 
  u_{m}<p_{n+1}\}|\ \mbox{ for } x\in [p_{n},p_{n+1}).$$
Then $g$ takes on 
decreasing values on $U\cap [p_{n},p_{n+1})$,  and $g\le f$,  because 
$|U\cap [p_{n},p_{n+1})|\le n<f(p_{n})$. Let $V\in\U$ be a set on which $g$ is 
nondecreasing: clearly $U\cap V$ has at most one point in each 
interval $[p_{n},p_{n+1})$.
 
\qed

\medskip    
    
    Following \cite{bdf}, let us consider the following families of 
    functions
	$$\S^{\U}=\{f \mid  \exists   U\in\U\ 
	\mbox{s.t.}\ f\, \mbox{1-1 on}\  U\,\},\
	\ \F_{\U}=\{f\,\mid\, \U\ \mbox{is}\ 
f\mbox{-QS}\,\},\ \ \mbox{and}$$
$$\G_{\U}=\{g\,\mid\, \exists U\in\U\,
\forall x,y\in U\ (x<y\ \; \Imp\ \; g(x)<y\,)\}.$$
Recall  the following facts, that represent three important features of QS 
ultrafilters, extensively used in \cite{bdf}:

\medskip\noindent
\textbf{Fact~1.}~(\cite[Theorem~1.1]{bdf})\ \emph{~If $\,\U$ is QS,
then\ $\F_{\U}=\G_{\U}$.}

\medskip\noindent
\textbf{Fact~2.}~(\cite[Proposition~1.5]{bdf}) \emph{ Let $g$ be 
interval-to-one, and put $g^{+}(x)= \max\,\{y\mid 
g(y)=g(x)\}$. Then $g\in\S^{\U}$ if and only if $g^{+}\in\G_{\U}$.
}

\medskip\noindent
\textbf{Fact~3.}~(\cite[Propositions~1.4~and~1.7]{bdf})  
\emph{$\F_{\U}$ is closed under sums, products, powers and compositions. 
Moreover $\G_{\U}$ has uncountable cofinality.
}

 \medskip
For general PQS ultrafilters we can prove both Facts 2 and 3, but 
only one
half of Fact 1, namely:
\begin{prop}\label{fac} ${}$Let $\,\U$ be PQS. Then
   \begin{enumerate}
\item 
\ $\F_{\U}\incl \G_{\U}$, and equality holds if and only if $\U$ is 
QS.

    \item For $g$ finite-to-one, put 
    $g^{+}(x)= \max\,\{y\mid\, g(y)=g(x)\}$: then $g\in\S^{\U}$
    if and only if $g^{+}\in\G_{\U}$.
    
    \item $\F_{\U}$ is closed under sums, products, powers and compositions; 
moreover $\G_{\U}$ has uncountable cofinality.

\end{enumerate}

\end{prop}

\proof ${}$

\noindent
1.  Assume that $\U$ is $f$-QS, with nondecreasing $f$, and pick any
   sequence $\la x_{n}\mid n\in\N\ra$ s.t.
   $x_{n+1}=f(x_{n})+x_{n}$.  Define the function  $h$ by 
   $h(x_{n}+j)=f(x_{n})-j$ for $0\le j < f(x_{n})$.
 Then there is a set in $A\in\U$ which meets each interval 
  $[x_{n},x_{n+1})$ in one point $a_{n}$. So by putting either 
  $u_{n}=a_{2n}$ or $u_{n}=a_{2n+1}$ we obtain a set $U\in\U$ 
  witnessing
  that $g=id+f$ belongs to $\G_{\U}$. Namely, in the even case we have
  $$u_{n+1}-u_{n}>x_{2n+2}-x_{2n+1}=f(x_{2n+1}) \ge f(u_{n}),$$
  and similarly in the odd case.

  The equality $\F_{\U}= \G_{\U}$ has been proved for QS 
  ultrafilters in Theorem 1.1 of \cite{bdf}. Finally,
the function $g$  has be choosen greater than the identity, so 
if $\U$ is not QS, then $g\notin\F_{\U}$, and the inclusion is \emph{proper}.

\medskip
\noindent
2.  Observe first that $g^{+}$ depends only on the partition induced 
by $g$, and not on its actual values. Moreover, if $h$ is  
any interval-to-one function inducing a coarser partition than $g$, then
$h^{+}\ge g^{+}$. Hence we may assume w.l.o.g. that $g$ is interval-to-one.

Assume $g^{+}\in\G_{\U}$, and pick $U=\{u_{n}\mid n\in\N\}\in\U$ 
such that
      $u_{n+1}>g^{+}(u_{n})$.  Suppose that $g(u_{n})=g(u_{n+1})$ 
       for some $n$: then 
       $g^{+}(u_{n})\ge u_{n+1} > g^{+}(u_{n})$, a
      contradiction. Hence $g$ is one-to-one on $U$.
      
      The reverse implication follows from the fact that 
      $g^{++}=g^{+}$.
      
      \medskip
      
      \noindent 
3. We prove first that if every function $g<f$ is 
$\U$-equivalent to a nondecreasing one, 
then the same property holds for every function $g<f^{2}$. 

Given $g$, let $h$ be the integral part of the square root of $g$. So 
$g<h^{2}+2h+1$, hence $g=h^{2}+h_{1}+h_{2}$ for suitable functions
$h_{1},h_{2}\le h<f$. By hypothesis we can pick nondecreasing 
functions $h',h'_{1},h'_{2}$ that are $\U$-equivalent to 
$h,h_{1},h_{2}$, respectively.
Then clearly $g$ is $\U$-equivalent to the nondecreasing function
$h'^{2}+h'_{1}+h'_{2}$.
So $\F_{\U}$ is closed under squares, and hence also under sums, 
products and powers. To settle compositions, observe first that, 
if $g, h\le id$, then $g\circ h\le h$, and the thesis is 
trivial. On the other hand, if $id\in\F_{\U}$, then $\U$ is QS, and 
we refer to the proof of Fact 3. given \emph{sub} Proposition 1.5 of
\cite{bdf}.

Finally, the proof of $\,\cof\G_{\U}>\og$ given \emph{sub} 
Proposition 1.7 of \cite{bdf} grounds 
solely on the fact that $\U$ is a P-point, so it works here as well.
\qed

\medskip
\emph{CAVEAT}:  
When $\U$ is not QS,  we may not state point 2  for
$\F_{\U}$, as it is done in \cite{bdf}, because $\G_{\U}$ is greater 
than $\F_{\U}$.

\medskip
The main tool in the study of PQS ultrafilters (and especially of PWR 
    ultrafilters in the next section) is the relative position of 
    particular cuts in the corresponding ultrapowers of $\N$.
    
Given a non-Q-point ultrafilter $\U$, let 
$\P=\la [p_{n},p_{n+1}) \mid n\in \N\ra$ be an interval partition 
witenssing the non-Q-pointness of $\U$, \ie\ such that there is  no $U\in\U$ 
with $|U\cap [p_{n},p_{n+1})|\le 1$ for all $n\in\N$.

For $U\in\U$ and $p_{n}\le x<p_{n+1}$ 
define the functions $a_{p}^{U},b_{p}^{U},$ and $c_{p}^{U}$ by
	$$a_{p}^{U}(x)=|U\cap [p_{n},x)|,\ \ b_{p}^{U}(x)=|U\cap 
     [x,p_{n+1})|,\ \ c_{p}^{U}(x)=|U\cap [p_{n},p_{n+1})|,$$
     and consider the corresponding families of functions
     $$\A_{p}^{\U} =\{a^{U}_{p} \mid    U\in\U \},\ \
     \B_{p}^{\U} =\{b^{U}_{p} \mid    U\in\U \},\ \ 
	\C_{p}^{\U} =\{c^{U}_{p} \mid    U\in\U \}.$$
	Put 
	$ \E^{\U}=\{f \mid  f\ \mbox{increasing}\!\!\! \mod\U
	\}$, and recall that 
	$\S^{\U}=\{f\,\mid\, f\ 1\mbox{-}1\!\!\! \mod\U\}.$

     %
%

 %
 %
\smallskip	
 We have

\begin{thm}\label{PQScut}
    Let $\,\U$ be a \emph{PQS} ultrafilter, and let $\P$ be
   an interval partition without selection set in $\U$. Let 
    $ F_{\U}$ be the cut
of the ultrapower 
$\N\ult{\N}{\U}$ whose left part is generated by $\F_{\U}$; 
let $E^{\U},S^{\U},A_{p}^{\U},B_{p}^{\U},$ and 
    $C_{p}^{\U}$ be the cuts
of the ultrapower 
$\N\ult{\N}{\U}$ whose right parts are generated by $\,\E^{\U},\S^{\U},
\A_{p}^{\U},\B_{p}^{\U},$ and 
$\C_{p}^{\U}$ respectively.

Then all cuts, but possibly $A_{p}^{\U}$, are greater than $\N$, and 
    $$A_{p}^{\U},S^{\U}\le E^{\U},\ \
    F_{\U}\le B_{p}^{\U},\ \  
\mbox{{and}}\ \    \max\{A_{p}^{\U},B_{p}^{\U}\}=C_{p}^{\U}.$$
 
  Moreover $\U$ is SPS
   if and only if $E^{\U}<F_{\U}$, and in 
  this case
$$ A_{p}^{\U}=S^{\U}= E^{\U} < F_{\U}\le  B_{p}^{\U}=C_{p}^{\U}.$$

\end{thm}

\proof
For $U\in\U$  put
$e^{U}(x)=|U\cap [0,x)|$, so every function increasing on $U$ is not 
smaller than $e^{U}$. Hence the cut $\E^{\U}$ is generated also by 
the set $\{e^{U} \mid    U\in\U \}$. Since $a_{p}^{U}\le e^{U}$, one 
gets $A_{p}^{\U}\le E^{\U}$.  The inequality $S^{\U}\le E^{\U}$ is 
trivial, and $F_{\U}\le B_{p}^{\U}$ holds because every  $U\in\U$ 
intersects some interval $[p_{n},p_{n+1})$ in more than one point, and hence
 no function $b_{p}^{U}$ is nondecreasing modulo $\U$.

 Moreover, for all\ $U\in\U$,
 $$ a^{U}_{p},b^{U}_{p}\le c^{U}_{p}= 
  a^{U}_{p}+b^{U}_{p},\ \
  \mbox{ whence } \ \ \frac{1}{2}c^{U}_{p}\le \max\{a_{p}^{U},b_{p}^{U}\}\le 
  c_{p}^{U}.$$
  Hence 
  $\max\{A_{p}^{\U},B_{p}^{\U}\}=C_{p}^{\U}$, because
  for all $U\in\U$ there exists $V\in\U$ s.t. 
  $c^{V}_{p}\le\frac{1}{2}c^{U}_{p}$.
  
One has $\N<F_{\U},
    S_{\U}$ because $\U$ is PQS, so it cannot be rapid. It follows 
    that only $A_{p}^{\U}$ might possibly be equal to $\N$.
    \smallskip
 
Finally, if 
$E^{\U}<F_{\U}$, then obviously $\S^{\U}\cap\F_{\U}\ne\0$. Conversely, 
$f\in \S^{\U}\cap\F_{\U}$ implies $f\in\E^{\U}$, and hence 
$ A_{p}^{\U}\le S^{\U}= E^{\U} < F_{\U}\le  B_{p}^{\U}=C_{p}^{\U}.$
Moreover if $a^{U}_{p}\in \F_{\U}$, then it is nondecreasing on some 
$V\in\U$. It follows that $a^{U}_{p}$ becomes increasing by taking 
off at most one point from each interval $V\cap[p_{n},p_{n+1})$, and 
the resulting set $V'$ belongs to $\U$, too. So $a^{U}_{p}\in 
\E^{\U}$, and also $A_{p}^{\U}= E^{\U}$. 

\qed

We conclude this section by
extending
Proposition 1.9 of \cite{bdf} to arbitrary PQS ultrafilters, thus
obtaining that the class of $f$-QS ultrafilters can be closed under 
isomorphisms only in the trivial case when every P-point is selective.

\begin{proposition}\label{bla1}
Assume that the ultrafilter $\U$ is not a Q-point, and let $f$ be 
an arbitrary nondecreasing unbounded function. Then there exists 
an increasing function $\fg$ such that 
the ultrafilter $\fg\U\cong \U$  is not $f$-QS.    
\end{proposition}

\proof
Let 
$\P=\la [p_{n},p_{n+1}) \mid n\in \N\ra$ be an interval partition 
witenssing the non-Q-pointness of $\U$, \ie\ such that there is  no $U\in\U$ 
with $|U\cap [p_{n},p_{n+1})|\le 1$ for all $n\in\N$.

Pick a sequence $b_{n}$\  such that\ ~$f(b_{n})> 
p_{n+1}$ ~and  ~$b_{n+1}-b_{n}> p_{n+1}-p_{n}$.
Define the function $\fg$ by 
$$\fg(p_{n}+j) = b_{n}+ j\ \mbox{ for }\  0\le j< p_{n+1}-p_{n}.$$
So the points $\fg(p_{n})=b_{n}$ determine an interval partition that has no 
selection set in $\fg\U$. Moreover $f(b_{n})> p_{n+1}$, hence
 any function $g$ such that $$g(b_{n}+j) =p_{n+1}-j\ \ \mbox{for}\ \
0\le j< a_{n+1}-a_{n}$$ is positive and not 
greater than $f$ on $\fg[\N]$, but cannot be nondecreasing modulo 
$\fg\U$.

\qed


\section{Weakly Ramsey ultrafilters}\label{wru}
 
An interesting weakening of the Ramsey property of selective 
ultrafilters has been considered by A. Blass in \cite{bl74}:

\begin{definition}
   \emph{The ultrafilter $\U$ on $\N$ is \emph{weakly Ramsey} (shortly WR) 
    if for
    every finite colouring of $[\N]^{2}$ there is $U\in\U$ s.t. $[U]^{2}$ 
    has \emph{only two colours}.\\
    $\U$ is \emph{properly weakly Ramsey} (shortly PWR) if it is WR, but not 
selective.}
    
\end{definition}

Throughout this section we assume that $\U$ is a PWR ultrafilter, and 
that $\P=\la [p_{n},p_{n+1}) \mid n\in \N\ra$ is an interval partition 
witnessing the non-selectivity of $\U$, so  there is  no $U\in\U$ 
with $|U\cap [p_{n},p_{n+1})|\le 1$ for all $n\in\N$.

	The behaviour 
	of functions modulo a PWR ultrafilter $\U$ is subject to severe 
	constraints,
	which recall those given by selectivity; namely every 
	function $f$ is $\U$-equivalent either to a $1$-to-$1$ 
	function, or to a function that is constant on each interval 
	$[p_{n},p_{n+1})$, independently of the choice of the 
	interval partition $\P$. More precisely (see Theorem 5 of \cite{bl74}):
	
\begin{lemma}\label{cresc}
    Let $f:\N\to\N$ and the 
    interval partition $\P$ be given. Then there exists $U\in\U$ such that exactly 
    one of the following cases occurs:\\
    ${}\ \ (i)$ $f$ is constant on $U$;\\
    ${}\ (ii)$  $f$ is increasing on $U$;\\
    $(iii)$  $f(x)< f(y)$ whenever  $x,y\in U$ 
    and there is $n$ such that $x<p_{n}\le y$, and  $f$ is constant on 
    $U\cap [p_{n},p_{n+1})$ for all $n\in \N$;\\
    ${}\ (iv)$ $f(x)< f(y)$ whenever  $x,y\in U$ 
    and there is $n$ such that $x<p_{n}\le y$, and  $f$ is decreasing on 
    $U\cap [p_{n},p_{n+1})$ for all $n\in \N$.
    
\smallskip    \noindent
    In particular,  
    the ultrafilter $f\U$ is selective if and only if $f$ is 
    constant on each interval $[p_{n},p_{n+1})$, \ie\ of type $(iii)$.
\end{lemma}
	
\proof
Put
$p(x)=n$ if $x\in [p_{n},p_{n+1}),$ and
identify $[\N]^{2}$ with the set of pairs 
$\{(x,y)\in\N^{2}\mid x<y\}$. Define the $6$-colouring of $[\N]^{2}$
according to all possible 
combinations of $p(x)\leq p(y)$ and $f(x)\gtreqless f(y)$.

\smallskip
By the choice of the interval partition, any $2$-coloured set 
$[U]^{2}$ with $U\in\U$ 
must comprehend both pairs with $p(x)= p(y)$ and pairs with $p(x)< p(y)$.
Now, when all are paired with $f(x)=f(y)$, then case (i) occurs, 
whereas case (ii) occurs when all are paired with $f(x)<f(y)$; case 
(iii) and (iv) occur when
$p(x)< p(y)$ is paired with $f(x)<f(y)$ and
$p(x)=p(y)$ with either $f(x)=f(y)$, or $f(x)> f(y)$, respectively.
It is easily seen that no one of the remaining cases can occur. \Pes, 
pairing $p(x)=p(y)$ with $f(x)<f(y)$ and $p(x)< p(y)$  with 
$f(x)=f(y)$ yields a contradiction by taking $p(x)=p(y)<p(z)$, \etc.

\smallskip
All functions of type $(ii)$ and $(iv)$ are $1$-$1$ modulo $\U$, so 
$f\U$ is isomorphic to $\U$. On the other hand,
 if $f$ is constant on 
each interval, then $g\circ f$ is non decreasing modulo $\U$ 
for all $g$. Hence all functions are nondecreasing modulo $f\U$, 
which is therefore selective.

\qed

\medskip

In order to classify  the different types of PWR 
ultrafilters, we 
recall the notation of Section \ref{fqsu}. 
For $U\in\U$ and $p_{n}\le x<p_{n+1}$ 
let 
	$$a_{p}^{U}(x)=|U\cap [p_{n},x)|,\ \ b_{p}^{U}(x)=|U\cap 
     [x,p_{n+1})|,\ \ c_{p}^{U}(x)=|U\cap [p_{n},p_{n+1})|;$$
     $$\A_{p}^{\U} =\{a^{U}_{p} \mid    U\in\U \},\
     \B_{p}^{\U} =\{b^{U}_{p} \mid    U\in\U \},\ \ 
	\C_{p}^{\U} =\{c^{U}_{p} \mid    U\in\U \}; $$
	$$\Sk^{\U}=\{f\,\mid\, f\ 1\mbox{-}1\!\!\! \mod\U\},\ \ 
	\ \ \mbox{and}\ \
	\E^{\U}=\{f \mid   f\ \mbox{increasing}  \!\!\! \mod\U
	\}.$$
	Then we have

 %
 %

\begin{thm}\label{PWRcut}
    Let $\U$ be a \emph{PWR} ultrafilter. Let
   $A_{p}^{\U},B_{p}^{\U},C_{p}^{\U},E^{\U},$ and 
    $S^{\U}$ be the cuts
of the ultrapower 
$\N\ult{\N}{\U}$ whose right parts are generated by
$\,\A_{p}^{\U},\B_{p}^{\U},\C_{p}^{\U},\E^{\U},$ and 
$\Sk^{\U}$ respectively.
Let  $ F_{\U}$ be the cut
of the ultrapower 
$\N\ult{\N}{\U}$ whose left part is 
generated by $\F_{\U}=\{f\,\mid\,\U\ f\mbox{-QS}\}$.
    Then, independently of the chosen interval partition,
 %
 %
 %
   $$\min\{A_{p}^{\U},B_{p}^{\U}\}= S^{\U}\le
   \begin{cases} A_{p}^{\U}= E^{\U}\\
       B_{p}^{\U}=F_{\U} 
    \end{cases}
    \le C_{p}^{\U}=\max\{A_{p}^{\U},B_{p}^{\U}\}$$
    
    Moreover $\U$ is rapid if and only if $\N= F^{\U}$, and then
     all considered cuts coincide with $\N$.

\end{thm}

\proof
According to Lemma \ref{cresc}, all functions are nondecreasing 
modulo $\U$, but those of type $(iv)$. Moreover every function of type 
$(iv)$ w.r.t. $U\in\U$ is greater than  $b^{U}_{p}$, so the cuts 
$F_{\U}$ and $B^{\U}_{p}$ coincide.

Similarly a function is $1$-$1$ on some $U\in\U$ if and only if its 
type is
either
$(ii)$ or $(iv)$. All functions of the former type are not less than  
 the corresponding function $a^{U}_{p}$, while those of the latter 
 type
 are not less than  
 the corresponding function $b^{U}_{p}$. Hence the cut $S^{\U}$ 
 coincides with the smaller between $A^{\U}_{p}$ and $B^{\U}_{p}$.
 
 The equality 
 $\max\{A_{p}^{\U},B_{p}^{\U}\}=C_{p}^{\U}$ has  been proved 
 in Theorem \ref{PQScut}, without any use of quasi-selectivity,
 as well as the trivial inequality $A_{p}^{\U}\le 
 E^{\U}$. On the other hand,
 each function $a_{p}^{U}$ is increasing modulo $\U$, so 
 for all $U\in\U$ there is $V\in\U$ such that $a_{p}^{U}\ge e^{V}$ on 
 $V$, and the converse inequality $A_{p}^{\U}\ge 
 E^{\U}$ follows.
 
  Finally, $\U$ being a P-point, it is rapid if and only if the 
  functions that are $1$-$1$ modulo $\U$ are coinitial in 
  $\N\ult{\N}{\U}\/ \N$, \ie\ $\N= S^{\U}$. But then also $F_{\U}$ 
  has to be equal to $\N$, otherwise $\U$ would be $f$-QS for some 
  $f$, and so selective by Proposition \ref{rap}. So it remains
  to prove that $\N= F_{\U}$ implies $\N= C_{p}^{\U}$. Assume 
  the contrary: then  $C_{p}^{\U}=A_{p}^{\U}>B_{p}^{\U}=\N$. 
  Define the bijection $\sg$ of $\N$ by $$\sg(x)=p_{n}+ p_{n+1}-x-1\ 
  \mbox{for }\ p_{n}\le x<p_{n+1}.$$ 
  Then clearly
  $$ a^{U}_{p}>_{\U} b^{U}_{p}\ \ \Iff\ \  a^{\sg U}_{p}<_{\sg\U} 
  b^{\sg U}_{p}.$$
  So $A_{p}^{\sg\U}< 
  B_{p}^{\sg\U}=F_{\sg\U}$, and $\sg\U\cong \U$ would be 
  simultaneously rapid and PQS, against Proposition \ref{rap}.
 
\qed

\bigskip
It follows immediately that a PWR ultrafilter $\U$ is QS if and only 
if the identity is less than the cut $B^{\U}$. More generally,
the above theorem allows for a complete specification of the ``quasi-selectivity'' 
properties of PWR ultrafilters. Namely

\begin{cor}\label{QS}
 Let $\U$ be a \emph{PWR} ultrafilter, and let 
    $A_{p}^{\U},B_{p}^{\U},$ and 
    $C_{p}^{\U}$ be the cuts
of the ultrapower 
$\N\ult{\N}{\U}$ whose right parts are generated by $\A_{p}^{\U},\B_{p}^{\U},$ and 
$\C_{p}^{\U}$ respectively.
  Then 

\begin{enumerate}

    \item  $\U$ is PQS  if and only if $\N\ne
    B_{p}^{\U}$, or equivalently if and only if $\,\U$ is not rapid;

    \item  $\U$ is SQS if and only if $ A_{p}^{\U}<B_{p}^{\U}$, or 
    equivalently
    $ A_{p}^{\U}\ne C_{p}^{\U}$;\\ (in particular $\U$ is QS if and 
    only if $id< C_{p}^{\U}$)
    
 \item  $\,\U$ is isomorphic to a QS ultrafilter if and only if\, 
$A_{p}^{\U}\ne B_{p}^{\U}$.

\end{enumerate} 
          
\end{cor}

\proof${}$

\noindent
1. Any  unbounded function $f<B^{\U}=F_{\U}$ witnesses that $\U$ is $f$-QS, and 
the last assertion 
of Theorem \ref{PWRcut} implies that such a function $f$ exists
unless $\U$ is rapid.

\medskip
\noindent
2. We have $C_{p}^{\U}=\max\{A_{p}^{\U},B_{p}^{\U}\}$, hence  
$A_{p}^{\U}\ne C_{p}^{\U}$ is equivalent to $S^{\U}=A_{p}^{\U}< 
B_{p}^{\U}=F_{\U}$, by Theorem \ref{PWRcut}. So there is $U\in\U$ s.t. 
$a^{U}_{p}<B^{\U}_{p}$: then $a^{U}_{p}$ is increasing modulo $\U$, 
and  $\U$ is $a^{U}_{p}$-QS.

\medskip
\noindent
3. If $A_{p}^{\U}< 
B_{p}^{\U}$, then $\U$ is SQS; so there is a function $f$ increasing 
modulo $\U$ such that $\U$ is $f$-QS. 
Then $\U\cong f\U$, and $f\U$ is QS by Proposition \ref{fQS}.

If $A_{p}^{\U}>
B_{p}^{\U}$, define the bijection $\sg$ of $\N$ by $\sg(x)=p_{n}+ p_{n+1}-x-1$ 
for $p_{n}\le x<p_{n+1}$. 
Then clearly
$$ a^{U}_{p}>_{\U} b^{U}_{p}\ \ \Iff\ \  a^{\sg U}_{p}<_{\sg\U} 
b^{\sg U}_{p}.$$
So $A_{p}^{\sg\U}< 
B_{p}^{\sg\U}$, and $\sg\U$ is isomorphic to a QS ultrafilter by the 
preceeding case.

Conversely, let  $\fg$ be a $1$-$1$ function, which we may assume of 
type $(ii)$ or $(iv)$, according to Lemma \ref{cresc}. In both cases 
there is an 
 interval partition $\P'$ such 
that $\fg[p_{n},p_{n+1})\incl [p'_{n},p'_{n+1})$ for all $n\in\N$. 
Then one has
$$a^{U}_{p}>_{\U} b^{U}_{p}\ \ \Iff\ \  a^{\fg U}_{p'}<_{\fg\U} 
b^{\fg U}_{p'}, \ \ \mbox{when}\ \fg \ \mbox{is of type}\ (iv);$$
whereas
$$a^{U}_{p}<_{\U} b^{U}_{p}\ \ \Iff\ \  a^{\fg U}_{p'}<_{\fg\U} 
b^{\fg U}_{p'}, \ \ \mbox{when}\ \fg \ \mbox{is of type}\ (ii).$$

It follows that the equality $A_{p}^{\U}=
B_{p}^{\U}$ is preserved under isomorphism, and such ultrafilters cannot 
be QS (nor SQS).

\qed

\bigskip

\section{Final remarks and open questions}\label{froq}

Recall that
    both 
    PWR and PQS ultrafilters are nonselective P-points, so the above results are nontrivial 
    only when such 
    ultrafilters exist. (And their existence is independent of \zfc\ 
    by a celebrated result of Shelah's, see \pes\ \cite{wi}.) However mild hypotheses, like \CH\ or \MA, 
    suffice in making both classes
    rich and distinct (see \cite{bl74,bdf}). 
   In fact 
these  classes are already different unless both are empty, because the former is closed under isomorphism, whereas the 
latter is not, by Proposition \ref{bla1}.

In \zfc, one can draw the following diagram of implications
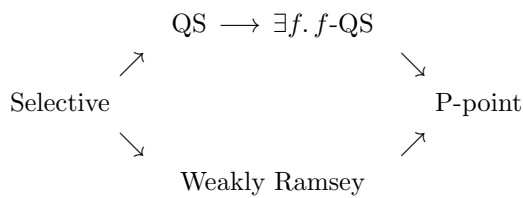
\begin{eqnarray*}
    {} & \mbox{QS}\,\longrightarrow\, \exists f.\,\mbox{$f$-QS} & {}  \\
    \nearrow & {} &\searrow \\
   \mbox{Selective}\ \ \ \   &{}  & \ \ \ \  \mbox{P-point}  \\
   \searrow & {} & \nearrow \\
    {} & \mbox{Weakly Ramsey} & {}
\end{eqnarray*}

 Recall that, assuming \CH, the following facts hold:
\begin{itemize}
    \item[(A)]  there exist  PWR ultrafilters $\U$ such that the cut induced by 
    $C_{p}^{\U}$ in the ultrapower $\N\ult{\N}{p\U}$ is arbitrarily 
    chosen among those having left part closed under exponentiation 
    and right part of uncountable coinitiality (Theorem $4$ of 
    \cite{bl74});\footnote{~It is worth mentioning that, according to 
    Theorem \ref{PWRcut}, if 
    $C_{p}^{\U}$ is taken to be 
    $\N$, then $\U$ is a \emph{rapid nonselective P-point}. Thus one has
     a non-forcing proof of the consistency of the existence 
    of such ultrafilters.}
    
    \item[(B)]  there are non-WR P-points (Theorem 2 of 
	\cite{bl74});
	
\item[(C)] 
there exist 
P-points that are not QS, and 
QS ultrafilters that are not selective (Theorem 1.2 of 
	\cite{bdf}).	
    
\end{itemize}

It follows from (A) that there exist rapid PWR ultrafilters, 
    necessarily not PQS, and also that for every $f$ there exist 
    $f$-QS PWR ultrafilters, necessarily non-$g$-QS for suitable 
    $g$.
    
 So, considering also (B-C), we may conclude that,   in the diagram above, 
 no arrow can be reversed nor inserted, except compositions.

Remark that both SQS and WR ultrafilters are  
P-points of a special kind, since they share the property that 
every function is equivalent to an
\emph{interval-to-one} function. So the question naturally arises as to 
whether this class of ``interval P-points'' is distinct from either 
one
of the other three classes. (We do not even know whether there exist 
WQS ultrafilters that are not ``interval P-points''.)

Many weaker conditions than the Continuum Hypothesis have been 
considered in the literature, in order to get more information 
about special classes of 
ultrafilters on $\N$. Of particular interest are (in)equalities among 
the so called ``combinatorial cardinal characteristics of the Continuum''. 
(\Pes\ one has that P-points or selective 
ultrafilters are generic if $\ck=\dk$ or $\ck=\cov(\B)$, respectively.
Moreover if $\cov(\B)<\dk=\ck$ then there are filters that are included 
in P-points, but cannot be 
extended to selective ultrafilters. See the 
comprehensive survey \cite{bl}.)
We conjecture that similar hypotheses can settle the problems 
mentioned above.

\bigskip

\bigskip
\bigskip

\end{document}